\documentclass[11pt]{article}      
\usepackage{graphicx}
\usepackage{latexsym}
\usepackage{amssymb}
\usepackage[mathscr]{eucal}
\usepackage{color}
\newcommand{\propref}[1]{{\rm Proposition \ref{proposition:#1}}}
\newcommand{\Bp}{\begin{proposition}\begin{sl}}
\newcommand{\EP}[1]{\end{sl}\label{proposition:#1}\end{proposition}}

\newcommand{\Div}{\mbox{\rm div}\,}

\newcommand{\supp}{\mbox{\rm supp}\,}

\newcommand{\curl}{\mbox{\rm curl}\,}

\newcommand{\Int}[2]{{\displaystyle \int_{ #1}^{ #2}}}

\newcommand{\Sup}[1]{{\displaystyle \sup_{#1}}}

\newcommand{\Frac}[2]{\displaystyle{\frac{\displaystyle{#1}}{\displaystyle{#2}}}}

\newcommand{\beea}{\begin{eqnarray}}

\newcommand{\eeea}{\end{eqnarray}}

\newcommand{\ms}{\medskip\smallskip}

\newcommand{\bfe}{{\mbox{\boldmath $e$}} }

\newcommand{\bfz}{{\mbox{\boldmath $z$}} }

\newcommand{\0}{{\mbox{\boldmath $0$}} }

\newcommand{\BF}{\begin{footnotesize}}

\newcommand{\EF}{\end{footnotesize}}

\newcommand{\pde}[2]{{\displaystyle \frac{\mbox{$\partial #1$}}{\mbox{$\partial #2$}}}}

\newcommand{\ode}[2]{{\displaystyle \frac{\mbox{$d #1$}}{\mbox{$d #2$}}}}

\newcommand{\bi}{\begin{itemize}}

\newcommand{\ei}{\end{itemize}}

\newcommand{\ed}{\end{document}}

\newcommand{\be}{\begin{equation}}

\newcommand{\ba}{\begin{array}}

\newcommand{\ea}{\end{array}}

\newcommand{\ee}{\end{equation}}

\newcommand{\eeq}[1]{\label{eq:#1}\end{equation}}

\newcommand{\real}{{\mathbb R}}

\newcommand{\nat}{{\mathbb N}}

\newcommand{\bfchi}{\mbox{\boldmath $\chi$}}

\newcommand{\bfomega}{\mbox{\boldmath $\omega$}}

\newcommand{\bfxi}{\mbox{\boldmath $\xi$}}

\newcommand{\bfx}{\mbox{\boldmath $x$}}

\newcommand{\bfy}{\mbox{\boldmath $y$}}

\newcommand{\bfphi}{\mbox{\boldmath $\varphi$}}

\newcommand{\bfv}{{\mbox{\boldmath $v$}} }

\newcommand{\bfu}{{\mbox{\boldmath $u$}} }

\newcommand{\bfw}{{\mbox{\boldmath $w$}} }

\newcommand{\bff}{{\mbox{\boldmath $f$}} }

\newcommand{\bfh}{{\mbox{\boldmath $h$}} }

\newcommand{\cala}{{\cal A}}

\newcommand{\bfF}{{\mbox{\boldmath $F$}} }

\newcommand{\bfb}{{\mbox{\boldmath $b$}} }

\newcommand{\bfg}{{\mbox{\boldmath $g$}} }

\newcommand{\bfn}{{\mbox{\boldmath $n$}} }

\newcommand{\half}{\mbox{$\frac{1}{2}$}}

\def\Bbb R{\real}

\def\tilde{\widetilde}

\def\bar{\overline}

\newcommand{\bfGamma}{\mbox{\boldmath $\Gamma$}}

\newcommand{\bfcalb}{\mbox{\boldmath ${\cal B}$}}

\newcommand{\bfcalg}{\mbox{\boldmath ${\cal G}$}}

\newcommand{\bfcalh}{\mbox{\boldmath ${\cal H}$}}

\newcommand{\bfcalf}{\mbox{\boldmath ${\cal F}$}}

\newcommand{\ED}{\end{description}}

\newcommand{\Footnote}{~\footnote}

\newcommand{\Br}{\begin{rem}\begin{rm}}

\newcommand{\Er}{\end{rm}\end{remark}}

\newtheorem{lemm}{Lemma}[section]

\newtheorem{theo}{Theorem}[section]

\newtheorem{rem}{Remark}[section]

\newtheorem{coro}{Corollary}[section]

\newtheorem{exe}{\footnotesize{Exercise}}[section]

\newcommand{\Be}{\begin{exe}\begin{footnotesize}\begin{rm}}

\newcommand{\EE}[1]{\end{rm}\end{footnotesize}\label{exe:#1}\end{exe}}

\newcommand{\Bt}{\begin{theo}\begin{sl}}

\newcommand{\Et}{\end{sl}\end{theorem}}

\newcommand{\Bl}{\begin{lemm}\begin{sl}}

\newcommand{\El}{\end{sl}\end{lemma}}

\newtheorem{proposition}{Proposition}[section]

\newcommand{\eqref}[1]{{\rm (\ref{eq:#1})}}

\newcommand{\Bc}{\begin{coro}\begin{sl}}

\newcommand{\Ec}{\end{sl}\end{coro}}

\newcommand{\ET}[1]{\end{sl}\label{theo:#1}\end{theo}}

\newcommand{\EL}[1]{\end{sl}\label{lemm:#1}\end{lemm}}

\newcommand{\theoref}[1]{{\rm Theorem \ref{theo:#1}}}

\newcommand{\ER}[1]{\end{rm}\label{rem:#1}\end{rem}}

\newcommand{\EC}[1]{\end{sl}\label{coro:#1}\end{coro}}

\newcommand{\lemmref}[1]{{\rm Lemma \ref{lemm:#1}}}

\newcommand{\essup}[1]{{\rm ess}\,{{\displaystyle \sup_{\hspace*{-6mm}{#1}}}}\,}

\begin{document}
\title{Existence, Uniqueness and Asymptotic Behavior\\ of Regular Time-Periodic Viscous Flow\\ around a Moving Body}
\author{Giovanni P Galdi \smallskip \\ { \small Department of Mechanical Engineering and Materials Science}\\ { \small University of Pittsburgh, USA}}
\date{}
\maketitle
\begin{abstract}We show existence and uniqueness of regular time-periodic solutions to the Navier-Stokes problem in the exterior of a rigid body, $\mathscr B$, that moves by arbitrary (sufficiently smooth) time-periodic translational motion of the same period, provided the size of the data is suitably restricted. Moreover, we characterize the spatial asymptotic behavior of such solutions and prove, in particular, that if $\mathscr B$ has a nonzero net motion identified by a constant velocity  $\bar{\bfxi}$ (say), then the solution exhibit a wake-like behavior in the direction $-\bar{\bfxi}$ entirely analogous to that of a  steady-state flow around a body that moves with  velocity $\bar{\bfxi}$.
\end{abstract}
\renewcommand{\theequation}{\arabic{section}.\arabic{equation}}
\setcounter{section}{0}
\section*{Introduction} 
Rigorous mathematical analysis of time-periodic flow of a Navier-Stokes liquid $\mathscr L$, around a {\em moving} rigid body, $\mathscr B$, is a relatively recent area of research.~\footnote{If the body is fixed, we refer the reader to \cite{Ma,MaPa,KoNa,Y,GaSo,KMT}.} In fact, the first contribution, due to A.L.~Silvestre and the present author, can be found in \cite{GS1}. In that paper the authors considered the general case where $\mathscr B$ moves by  arbitrary motion characterized by (sufficiently smooth)     time-periodic  translational velocity $\bfxi=\bfxi(t)$, and angular velocity $\bfomega=\bfomega(t)$. In particular, they  showed existence of corresponding solutions to the associated Navier-Stokes problem in a ``weak" class (a la Leray-Hopf) for data of arbitrary size, and in a ``strong" class (a la Ladyzhenskaya) if the size of the data is appropriately restricted. However, the important problem of uniqueness of these solutions was left open. \par
The question was successively reconsidered and thoroughly  investigated by a number of authors who, by entirely different methods, were able to prove existence and uniqueness of time-periodic solutions of period $T$ (from now on referred to as ``$T$-periodic" solutions) in appropriate function classes, under the assumption that both characteristic vectors $\bfxi$ and $\bfomega$ are constant \cite{NTH,GH,HG,GaKy,Ngu,GaKy1,MH}, and a $T$-periodic body force is acting on $\mathscr L$.

Very recently, in  \cite{GARMA} we began to investigate the above properties in the general situation when $\bfxi$ is not  constant, while assuming    $\bfomega\equiv\0$ . Converted in mathematical terms, this amounts to find $T$-periodic solutions $(\bfu,p)$ to the following system of equations
\be\ba{cc}\smallskip\left.\ba{ll}\medskip
{\partial}_t\bfu-\bfxi(t)\cdot\nabla\bfu+\bfu\cdot\nabla\bfu=\Delta\bfu-\nabla {p}+\bfb\\
\Div\bfu=0\ea\right\}\ \ \mbox{in $\Omega\times (-\infty,\infty)$}\\
\bfu(x,t)=\bfxi(t)\,,\ \ (x,t)\in \partial\Omega\times (-\infty,\infty)\,,
\ea
\eeq{0.1}
where $\bfu$ and $p$ are velocity and pressure fields of $\mathscr L$,   and $\Omega$ is the exterior of a connected compact region  of $\real^3$ (the body $\mathscr B$). Moreover, for the sake of generality, we include also a (prescribed) $T$-periodic body force $\bfb=\bfb(x,t)$ acting on $\mathscr L$. 
In \cite{GARMA} we analyzed the case when $\bfxi$ has zero average over a period, namely,  \be\bar{\bfxi}:=\frac1T\int_0^T\bfxi(t)\,{\rm d}t=\0,\eeq{0} which implies that  $\mathscr B$ has zero net motion. This happens, for example, if $\mathscr B$ oscillates between two fixed configurations. We then showed, in particular,   existence, uniqueness and regularity of such a flow, on condition that  $\bfxi$ is, in suitable norm, below a certain constant depending only on  $\Omega$ and $T$. Furthermore, we proved that  $\bfu(x,t)$ decays like $|x|^{-1}$, uniformly in time $t\in [0,T]$, where $|x|$ denotes the distance of a generic point in $\Omega$ from the origin located in  $\mathscr B$. Notice that this behavior is the same as that of a steady-state flow around an immovable body \cite[Section X.9]{Gab}.\par
Objective of the present work is to continue and, to an extent, complete the research initiated in \cite{GARMA}. Specifically, we shall investigate the same problem as in \cite{GARMA}, but relaxing the assumption \eqref{0}, thus  allowing $\mathscr B$ to have a nonzero net motion over a period of time. Also for this more general problem, we are able to  show  existence and uniqueness of regular solutions if the data are suitably restricted. Concerning the asymptotic spatial behavior, we demonstrate the following. Without loss of generality, take $\bar{\bfxi}=\lambda\,\bfe_1$, with $\bfe_1$ unit vector in the direction $x_1$ and $\lambda\ge 0$. Then $\bfu(x,t)$ decays like $|x|^{-1}[1+\lambda\,(|x|+x_1)]^{-1}$, uniformly in $t\in[0,T]$. As a consequence, if $\lambda=0$ we find a  result in agreement with that obtained in \cite{GARMA}. However, if $\lambda>0$, the velocity field shows a ``wake" behavior in the direction opposite to $\bar{\bfxi}$,  entirely analogous to that of a  steady-state flow around a body that moves with constant velocity $\bar{\bfxi}$  \cite[Section X.8]{Gab}. 
\par
The method we use here is a generalization of that employed in \cite{GARMA} and relies upon the proof of existence,  uniqueness and corresponding estimates of solutions to the linear counterpart of problem \eqref{0.1} in a specific function class. Members of this class are regular in a well-defined sense, on the one hand, and, on the other hand, they decay at large spatial distances uniformly in time in a suitable fashion, provided the data decay appropriately as well; see \propref{1}.  With this result in hand, it is then quite straightforward to apply the contraction mapping theorem and prove analogous results for the full nonlinear problem \eqref{0.1} on condition that $\bfxi$ and $\bfb$ are below a constant depending on $\Omega$ and $T$; see \theoref{3.1}.  
\par
The plan of the paper is as follows. After recalling some preliminary lemmas in Section 1, in the following Section 2 we prove the well-posedness results mentioned above for the linear problem obtained from \eqref{0.1} by neglecting the nonlinear term $\bfu\cdot\nabla\bfu$ and replacing $\bfb$ with a function $\bff$ satisfying suitable regularity and spatial decay conditions. This result is achieved in two steps. In the first one, we construct unique regular solutions by combining Galerkin's method with the classical ``invading domains" procedure; see \lemmref{4.2}. This finding requires that $\bfxi$ and $\bff$ possess a certain degree of smoothness. If, in addition, $\bff$ decays at large distances and uniformly in time at a suitable rate, we then show that the above solutions must decay as well;  see \lemmref{1}. The two lemmas are then combined in \propref{1}, to provide the desired well-posedness result. In the final Section 3, we employ \propref{1} in combination with a classical fixed-point argument to extend the results of that proposition to the full nonlinear problem \eqref{0.1}, under suitable restrictions on the magnitude of $\bfxi$ and $\bfb$ in suitable norms; see \theoref{3.1}.    
\section{Preliminaries}
We begin to recall some  notation.
Throughout, $\Omega$ denotes  the complement of the closure of a bounded  domain $\Omega_0\subset\mathbb R^3$, which we  assume   of class $C^2$, and  take the origin of the coordinate system in the interior of $\Omega_0$.   For $R\ge R_*:=2{\rm diam}\,(\Omega_0)$, we set 
$
\Omega_R=\Omega\cap \{|x|<R\}\,,\ \ \Omega^R=\Omega\cap\{|x|>R\}$.
For a given domain $A\subseteq \real^3$, by $L^q (A)$,  $1\leq q \leq \infty,$  
$W^{m,q}({A}),$ $W_0^{m,q}(A)$, $m \geq 0,$  $(W^{0,q}\equiv W^{0,q}_0\equiv L^q$), we denote usual Lebesgue and Sobolev classes, with corresponding norms $\|.\|_{q,A}$ and $\|.\|_{m,q,A}$.\Footnote{We shall use the same font style to denote scalar, vector and tensor
function spaces.} The letter $P$ stands for the (Helmholtz) projector from $L^2(A)$ onto its subspace constituted by solenoidal (vector) function with vanishing normal component, in distributional sense, at $\partial A$.   
 We also set $\int_{A}u\cdot v=\langle u,v\rangle_{A}$.  
$D^{m,2}(A)$ is the space of (equivalence classes of) functions $u$ such that
$ 
\sum_{|k|=m}\|D^k u\|_{2,A}<\infty\,.
$ 
Obviously, the latter defines a seminorm in $D^{m,2}(A)$. Also,  by $D_0^{1,2}(A)$ we denote the completion of $C_0^\infty(A)$ in the norm $\|\nabla(\cdot)\|_2$. In the above notation,  the subscript ``$A$" will be omitted, unless confusion arises. A function $u:A\times \real\mapsto \real^3$ is 
{\em $T$-periodic}, $T>0$, if $u(\cdot,t+T)=u(\cdot\,t)$, for a.a. $t\in \real$,
 and we set
$
{\bar u}:=\frac{1}{T}\int_{0}^{T}u(t){\rm d}t\,.
$
Let $B$ be a function space endowed with seminorm $\|\cdot\|_B$, $r=[1,\infty]$, and $T>0$. $L^r(0,T;B)$ is the class of functions
$u:(0,T)\rightarrow B$ such that 
$$
\|u\|_{L^r(B)}\equiv\left\{\ba{ll}\smallskip\big( \Int{0}{T}\|u(t)\|_B^r \big)^{\frac 1r}<\infty, \ \ \mbox{if 
$r\in [1,\infty)\,;$}\\   
\essup{t\in[0,T]}\,\|u(t)\|_B <\infty, \ \ \mbox{if $r=\infty.$}
\ea\right.
$$
Likewise, we put
$$
W^{m,r}(0,T;B)=\Big\{u\in L^{r}(0,T;B): \textcolor{black}{\partial_t^ku\in L^{r}(0,T;B), \, k=1,\ldots,m}\Big\}\,.
$$
Unless confusion arises, we shall simply write $L^r(B)$ for $L^r(0,T;B)$, etc. 
Finally, if $A:=\Omega,\real^3$,  $m\ge 1$, and $\lambda\ge 0$ we set
$$\ba{ll}\medskip
[\!]f[\!]_{m,\lambda,A}:=\Sup{x\in A}\,|(1+|x|)^{m}(1+2\lambda \,s(x))^m f(x)|\,,\\ 

[\!] f [\!]_{\infty,m,\lambda,A}:=\Sup{(x,t)\in A\times (0,\infty)}\,|(1+|x|)^{m}(1+2\lambda \,s(x))^m f(x,t)|\,.\ea
$$
where $s(x)=|x|+x_1$,  $x\in\real^3$, and  the subscript $A$ will be omitted, unless necessary.
\par
We next collect some preliminary results whose proof is given elsewhere. We begin with the following one,  a special case of \cite[Lemma II.6.4]{Gab}
\Bl
There exists a  function $\psi_{R} \in C_0^{\infty}(\real^n)$ defined for all $R>0$ such that $0\le \psi_{R}(x)\le 1$, $x\in \real^n$, and satisfying the following properties 
$$
\displaystyle \lim_{R\to\infty}\psi_{R}(x)=1\,,\  \mbox{uniformly pointwise\,;}\ \ \
\left\|\pde{\psi_{R}}{x_1}\right\|_{\frac32}\le C_1\,,
$$
where $C_1$  is independent of  $R.$ Moreover, the support of $\partial\psi_{R}/\partial x_j$, $j=1,\ldots,n$, is contained in $\Omega^{\frac{R}{\sqrt2}}$ and
$$
\|u\,|\nabla\psi_{R}|\,\|_{2}\le C_2\,\|\nabla u\|_{2,
\Omega^{\frac{R}{\sqrt{2}}}}\,,\ \ \mbox{for all $u\in D_0^{1,2}(\Omega)$}\,.
$$
where  $C_2$ is independent of $R$. 
\EL{1.1}
The following result can be found  in \cite[Exercise III.3.7]{Gab}.
\Bl Let $\cala$ be a bounded Lipschitz domain in $\real^3$, and let $f\in W_0^{1,2}(\cala)$ with $\int_{\cala}f=0$. Then the problem
\be
\Div\bfz=f\,\ \mbox{in $\cala$}\,, \ \bfz\in W_0^{2,2}(\cala)\,,\ \ \|\bfz\|_{2,2}\le C_0\,\|f\|_{1,2}\,,
\eeq{Bog}
for some $C_0=C_0(\cala)>0$ has at least one solution. Moreover, if $f=f(t)$ with $\partial_t f\in L^\infty(L^2(\cala))$, then we have also $\partial_t\bfz\in L^\infty(W_0^{1,2}(\cala))$ and
$$
\|\partial_t\bfz\|_{1,2}\le C_0\,\|\partial_t f\|_2\,.
$$ 
\EL{1.1_1}
The next result is proved in  \cite[Lemma 2.2]{GS1}. 
\Bl
Let $\bfxi \in W^{\textcolor{black}{2},2}(0,T)$ be $T$-periodic. There exists a solenoidal, $T$-periodic function $\tilde{\bfu} \in W^{\textcolor{black}{1},2}(W^{m,q}),$ $m\in \mathbb{N},$ $q\in[1,\infty],$ such that 
$$\ba{ll}\medskip
\tilde{\bfu}(x,t)=\bfxi(t)\,,\ (t,\bfx)\in[0,T]\times\partial\Omega\,,\\ \medskip
\tilde{\bfu}(x,t)=0\,,\ \mbox{for all $t\in [0,T]$, all $|\bfx|\ge\rho$, and some $\rho>R_*$}\,,\\ 
\| \tilde{\bfu} \|_{W^{2,2}(W^{m,q})}\leq C\,\|\bfxi\|_{W^{2,2}(0,T)}
\,,\ea
$$
where $C=C(\Omega,m,q)$. 
\EL{ext}
We conclude by recalling the following lemma  showing suitable existence and uniqueness properties for a  linear Cauchy problem \cite[Theorem {VIII.4.4}]{Gab}
\Bl
Let $\bfcalg$ be a second-order tensor field in $\mathbb{R}^3 \times (0,\infty)$ such that
$$
[\!]\bfcalg(t)[\!]_{\infty,2,\lambda} + \essup{t\geq 0}\|\nabla \cdot \bfcalg(t)\|_2 < \infty\,,
$$
and let $\bfh\in L^{\infty,q}(\real^3\times (0,\infty))$, $q\in (3,\infty)$, with spatial support contained in a ball of radius $\rho$, some $\rho>0$, centered at the origin.
Then, the problem 
\be
\ba{ll}\ms\left.\ba{rl}\ms
\pde{\bfw}{t}&=\Delta\bfw+\lambda\,\pde{\bfw}{x_1}-\nabla\phi+\nabla\cdot\bfcalg+\bfh\\
\nabla\cdot\bfw&=0\ea\right\}\ \ \mbox{in $\real^3\times (0,T)$}\\
\bfw(x,0)=\0\,,
\ea
\eeq{VIII.4.1}
has one and only one solution such that for all $T>0$,
\begin{equation}
\bfw\in L^2(0,T;W^{2,2})\,, \ \partial_t\bfw\in L^2(0,T;L^2)\,;\ \ \nabla\phi\in L^{2}(0,T;L^2).
\label{estimunn1}
\end{equation}
Moreover, 
$$
[\!] \bfw(t) [\!]_{\infty,1,\lambda}
<\infty\,,
$$
and the following inequality holds:
\begin{equation}
[\!] \bfw(t) [\!]_{\infty,1,\lambda} 
\leq C\,\left([\!]\bfcalg(t)[\!]_{\infty,2,\lambda} +\essup{t\geq 0}\|\bfh(t)\|_q\right)\, \,  
\label{estimunn2}
\end{equation}
with $C=C(q,\rho,B)$, whenever ${\lambda}\in [0,B]$, for some $B>0$. 
\EL{VIII.4.4}
\setcounter{equation}{0}
\section{On the Unique Solvability of the Linear Problem}
The  main objective of this section is to prove existence and uniqueness of $T$-periodic solutions, in appropriate function classes, to the following set of linear equations:
\be\ba{cc}\smallskip\left.\ba{ll}\medskip
{\partial}_t\bfu-\bfxi(t)\cdot\nabla\bfu=\Delta\bfu-\nabla {p}+\bff\\
\Div\bfu=0\ea\right\}\ \ \mbox{in $\Omega\times (0,T)$}\\
\bfu(x,t)=\bfxi(t)\,,\ \ (x,t)\in \partial\Omega\times [0,T]\,,
\ea
\eeq{4.7}
where $\bfxi=\bfxi(t)$ and $\bff=\bff(x,t)$ are suitably prescribed $T$-periodic functions. Without loss, we take
$$
\bar{\bfxi}=\lambda\,\bfe_1\,,\ \ \lambda\ge 0\,,
$$
where $\bfe_1$ is the unit vector along the $x_1$-axis.
\medskip\par We begin show the following result. 
\Bl Let 
$$\bff=\Div \bfcalf\in W^{1,2}(L^{2}(\Omega))\,, \ \mbox{$\bfcalf\in L^{2}(L^{2}(\Omega))$\,,}
$$
and $\bfxi\in W^{2,2}(0,T)$ be prescribed $T$-periodic functions. Then, there exists one and only one  $T$-periodic solution $(\bfu,p)$ to \eqref{4.7} such that
\be
\bfu\in L^\infty(L^6)\,,\  \partial_t\bfu\in  L^{\infty}(W^{1,2})\cap L^2(D^{2,2})\,, \ \nabla\bfu\in L^\infty(W^{1,2})\,;\ \ p\in L^\infty(D^{1,2})\,.
\eeq{class}
Furthermore, 
\be\ba{ll}\medskip
\|\partial_t\bfu\|_{L^{\infty}(W^{1,2})\cap L^2(D^{2,2})}+\|\bfu\|_{L^{\infty}(L^6)}+\|\nabla\bfu\|_{L^\infty(W^{1,2})} +\|\nabla p\|_{L^\infty(D^{1,2})}+\|p\|_{L^\infty(L^2(\Omega_R))}\\ 
\hspace*{1cm}\le C\,\big(\|\bff\|_{W^{1,2}(L^2)}+\|\bfcalf\|_{L^{2}(L^{2})}+\|\bfxi\|_{W^{2,2}(0,T)} \big)
\ea
\eeq{est}
where $C=C(\Omega,T,R,\xi_0)$, 
for any fixed $\xi_0$ such that $\|\bfxi\|_{W^{2,2}(0,T)}\le \xi_0$. 

\EL{4.2}
{\em Proof.} We follow the argument of \cite[Sections 3 \& 4]{GS1}  to show existence, by combining the classical Galerkin method with the ``invading domains" procedure. We shall  limit ourselves to prove the basic {\em a priori} estimates, referring the reader to that article for the (classical) procedure of how these estimates can be used to prove the stated existence result. Let $\bfu=\bfv+\tilde{\bfu}$, with $\tilde{\bfu}$ given in \lemmref{ext}, and  consider problem \eqref{4.7} along an increasing, unbounded sequence of (bounded) domains $\{\Omega_{R_k}\}$ with $\cup_{k\in\nat}\Omega_{R_k}=\Omega$, that is,
\be\ba{cc}\smallskip\left.\ba{ll}\medskip
{\partial}_t\bfv_k-\bfxi(t)\cdot\nabla\bfv_k=\Delta\bfv_k-\nabla \tilde{p}_k+\bff+\bff_c\\
\Div\bfv_k=0\ea\right\}\ \ \mbox{in $\Omega_{R_k}\times (0,T)$}\\
\bfv_k(x,t)=\0\,,\ \ (x,t)\in \partial\Omega_{R_k}\times [0,T]\,,
\ea
\eeq{4.8}
where 
$$
\bff_c:= \Delta\tilde{\bfu}-\partial_t\tilde{\bfu}+\bfxi(t)\cdot\nabla\tilde{\bfu}
$$
If we formally dot-multiply \eqref{4.8}$_1$ by $\bfv_k$ and integrate by parts over $\Omega_{R_k}$ we get
\be
\half\ode{}t\|\bfv_k(t)\|_2^2+\|\nabla\bfv_k(t)\|_2^2=\langle \bff+\bff_c,\bfv_k\rangle\le c_0\left(\|\bfcalf\|_2+\|\bff_c\|_{\frac65}\right)\|\nabla\bfv_k\|_2\,, 
\eeq{4.9}
where we have used  the Sobolev inequality 
\be
\|\bfz\|_6\le \gamma_0\,\|\nabla\bfz\|_2,\ \ \ \bfz\in D_0^{1,2}(\real^3)\,, 
\eeq{So}
with $\gamma_0$ numerical constant. 
Employing in \eqref{4.9} Cauchy inequality along with Poincar\`e inequality $\|\bfv_k\|_2\le c_{R_k}\|\nabla\bfv_k\|_2$  we get, in particular,
$$
\ode{}t\|\bfv_k(t)\|_2^2+c_{1R_k}\|\bfv_k(t)\|_2^2\le c_2\,\left(\|\bfcalf\|_2^2+\|\bff_c\|_{\frac65}^2\right)\,.
$$
Proceeding as in  \cite[Lemma 3.1]{GS1}, we may combine this inequality with Galerkin method to prove the existence of a $T$-periodic (distributional) solution $\bfv_k$ to \eqref{4.8} with $\bfv_k\in L^{\infty}(L^2(\Omega_{R_k}))\cap L^2(D_0^{1,2}(\Omega_{R_k}))$ .  In addition,
\be
\|\nabla\bfv_k\|_{L^{2}(L^{2})}\le c\,\left(\|\bfcalf\|_{L^{2}(L^{2})}+\|\bff_c\|_{L^{2}(L^{\frac65})}\right)\,,
\eeq{4.10}
where the constant $c$ is independent of $R_k$; see \cite[Section 3]{GS1} for  details. 
We will next show uniform (in $k$)  estimates for $\bfv_k$ in spaces of higher regularity. In this regard, we notice that 
by the mean value theorem,   from \eqref{4.10}  it follows that there is $t_0\in (0,T)$ such that
\be
\|\nabla\bfv_k(t_0)\|_2^2\le \,c_3\left(\|\bfcalf\|_{L^{2}(L^{2})}^2+\|\bff_c\|^2_{L^{2}(L^{\frac65})}\right)\,.
\eeq{mvt}
If we formally dot-multiply both sides of \eqref{4.8}$_1$ a first time by $P\Delta\bfv_k$, a second time by $\partial_t\bfv_k$ and integrate by parts over $\Omega_{R_k}$, we deduce
\be\ba{ll}\medskip
\half\ode{}t\|\nabla\bfv_k(t)\|_2^2+\|P\Delta\bfv_k(t)\|_2^2=\langle\bfxi\cdot\nabla\bfv_k,P\Delta\bfv_k\rangle+\langle\bff+\bff_c,P\Delta\bfv_k(t)\rangle\\
\half\ode{}t\|\nabla\bfv_k(t)\|_2^2+\|\partial_t\bfv_k(t)\|_2^2=\langle\bfxi\cdot\nabla\bfv_k,\partial_t\bfv_k\rangle+\langle\bff+\bff_c,\partial_t\bfv_k(t)\rangle\,.
 
\ea
\eeq{JB}
Therefore, summing side-by-side the two equations in \eqref{JB} and employing Cauchy-Schwarz inequality allows us to infer
\be
\ode{}t\|\nabla\bfv_k(t)\|_2^2+c_4(\big(\|\partial_t\bfv_k(t)\|_2^2+\|P\Delta\bfv_k(t)\|_2^2\big)\le c_4\big(\|\bff\|_2^2+\|\bff_c\|_2^2+\|\nabla\bfv_k(t)\|_2^2\big)\,, 
\eeq{KL}
with $c_4=c_4(\xi_0)$.
We now recall the inequality
\be
\|D^2\bfz\|_{2,\Omega_R}\le c_{\Omega}\,\left(\|P\Delta\bfz\|_{2,\Omega_R}+\|\nabla\bfz\|_{2,\Omega_R}\right)\,,\ \ \bfz\in D^{1,2}(\Omega_R)\cap D^{2,2}(\Omega_R)\,, 
\eeq{hey}
with $c_{\Omega}$ depending only on the regularity of $\Omega$ but {\em not} on $R$  \cite[Lemma 1]{Hey}. Thus,
integrating both sides of \eqref{KL}  over $[t_0,t]$,   using the $T$-periodicity property along with \eqref{mvt},  \eqref{hey}, and \lemmref{ext} we  show  that  $\bfv_k\in W^{1,2}(L^2(\Omega_{R_k}))\cap L^\infty(D_0^{1,2}(\Omega_{R_k}))\cap L^2(D^{2,2}(\Omega_{R_k}))$ and, in addition,  $\bfv_k$ satisfies the uniform bound \cite[Lemma 4.1]{GS1} 
\be\ba{rl}\medskip
\|\bfv_k\|_{L^{\infty}(L^6)}+\|\nabla\bfv_k\|_{L^{\infty}(L^2)}+&\|\partial_t\bfv_k\|_{L^{2}(L^2)}+\|D^2\bfv_k\|_{L^{2}(L^2)}\\ \medskip &\le c\,\big(\|\bff\|_{L^{2}(L^2)}+\|\bfcalf\|_{L^{2}(L^2)}+\|\bff_c\|_{L^{2}(L^\frac65)}\big)\\ &
\le C\,\big(\|\bff\|_{L^{2}(L^2)}+\|\bfcalf\|_{L^{2}(L^{2})}+\|\bfxi\|_{W^{2,2}(0,T)} \big)
\,,\ea
\eeq{4.11}
with $C$ independent of $R_k$. Next, we take the time derivative of both sides of \eqref{4.8}$_1$, and dot-multiply both sides of the resulting equation a first time by $\partial_t\bfv_k$, a second time by $P\Delta\partial_t\bfv_k$ and then integrate  over $\Omega_{R_k}$. We then obtain
\be
\half\ode{}t\|\partial_t\bfv_k(t)\|_2^2+\|\nabla\partial_t\bfv_k(t)\|_2^2=\langle\bfxi'\cdot\nabla\bfv_k,\partial_t\bfv_k\rangle+\langle\partial_t\bff+\partial_t\bff_c,\partial_t\bfv_k(t)\rangle\,,
\eeq{int}
and
\be\ba{rl}\medskip
\half\ode{}t\|\nabla\partial_t\bfv_k(t)\|_2^2+\|P&\!\!\!\!\!\Delta\partial_t\bfv_k(t)\|_2^2\\&=\langle\bfxi'\cdot\nabla\bfv_k,P\Delta\partial_t\bfv_k(t)\rangle+\langle\partial_t\bff+\partial_t\bff_c,P\Delta\partial_t\bfv_k(t)\rangle\,.
\ea\eeq{chi}
From \eqref{4.11} and the mean value theorem we find that there is  $t_1\in (0,T)$ such that
\be
\|\partial_t\bfv_k(t_1)\|_2\le c\,\big(\|\bff\|_{L^{2}(L^2)}+\|\bfcalf\|_{L^{2}(L^{2})}+\|\bfxi\|_{W^{2,2}(0,T)} \big)\,,
\eeq{spa}
and so, integrating \eqref{int} over $[t_1,t]$
and using Cauchy-Schwarz inequality, \eqref{spa},  \eqref{4.11} and the $T$-periodicity of $\bfv_k$, we arrive at
\be
\|\partial_t\bfv_k\|_{L^\infty(L^{2})}+\|\nabla\partial_t\bfv_k\|_{L^{2}(L^2)}\le C\,\big(\|\bff\|_{W^{1,2}(L^2)}+\|\bfcalf\|_{L^{2}(L^{2})}+\|\bfxi\|_{W^{2,2}(0,T)} \big)
\eeq{JB0}
By a similar token, from \eqref{chi},  \eqref{JB0} and \eqref{hey}, we get
\be\ba{rl}\medskip
\|\nabla\partial_t\bfv_k\|_{L^\infty(L^{2})}+\|&\!\!\!\!\!D^2\partial_t\bfv_k\|_{L^{2}(L^2)}\\ &\le C\,\big(\|\bff\|_{W^{1,2}(L^2)}+\|\bfcalf\|_{L^{2}(L^{2})}+\|\bfxi\|_{W^{2,2}(0,T)} \big)\,.
\ea\eeq{JB01}
Therefore, combining \eqref{4.11}, \eqref{JB0}, and \eqref{JB01} we infer
\be\ba{ll}\medskip
\|\partial_t\bfv_k\|_{L^\infty(W^{1,2})}+\|\bfv_k\|_{L^\infty(L^6)}+\|\nabla\bfv_k\|_{L^\infty(L^2)}+\|D^2\bfv_k\|_{W^{1,2}(L^2)}\\ \medskip
\hspace*{4cm}\le C\,\big(\|\bff\|_{W^{1,2}(L^2)}+\|\bfcalf\|_{L^{2}(L^{2})}+\|\bfxi\|_{W^{2,2}(0,T)} \big)
\ea
\eeq{JB1}
where $C$ is independent of $k$.
Finally, setting $\bfF_k:=\Delta\bfv_k+\bff+\bff_c$, from \eqref{4.8}$_1$ we get, formally, that $\tilde{p}_k$ obeys for a.a. $t\in [0,T]$ the following Neumann problem\footnote{Note that $\bfxi(t)\cdot\nabla\bfv_k\cdot\bfn|_{\partial\Omega_{R_k}}=0$.}
\be
\Delta \tilde{p}_k=\Div\bfF_k\ \ \mbox{in $\Omega_{R_k}$}\,,\ \ \partial \tilde{p}_k/\partial\bfn|_{\partial\Omega_{R_k}}=\bfF_k\cdot\bfn\,.
\eeq{4.12}  
Therefore, multiplying both sides of the first equation by $\tilde{p}_k$ and integrating by parts over $\Omega_{R_k}$ we easily establish that the pressure field ${p}_k$ associated to $\bfv_k$ satisfies the estimate \cite[Lemma 4.3]{GS1}
\be
\|\nabla \tilde{p}_k\|_2\le c\,\big(\|D^2\bfv_k\|_2+\|\bff\|_2+\|\bff_c\|_2\big)
\eeq{4.13}
with $c$ independent of $k$. We may now let $R_k\to\infty$ and use the uniform estimate \eqref{JB1} and \lemmref{ext}, 
to show the existence of a pair 
$
(\bfu:=\bfv+\tilde{\bfu},\tilde{p})$, with $\bfu$ $T$-periodic, in the class 
\be
\partial_t\bfu\in  L^{\infty}(W^{1,2})\cap L^{2}(D^{2,2})\,,\ \bfu\in L^{\infty}(L^{6})\,,\ \nabla\bfu\in L^{\infty}(L^{2})\cap L^2(D^{1,2})\,, \ \   \tilde{p}\in L^2(D^{1,2})\,,
\eeq{JB0}
such that
\be\ba{ll}\medskip
\|\partial_t\bfu\|_{L^\infty(W^{1,2})}+\|\bfu\|_{L^\infty(L^6)}+\|\nabla\bfu\|_{L^\infty(L^2)}+\|D^2\bfu\|_{W^{1,2}(L^2)}+\|\nabla\tilde{p}\|_{L^2(L^2)}\\ \medskip
\hspace*{4cm}\le C\,\big(\|\bff\|_{W^{1,2}(L^2)}+\|\bfcalf\|_{L^{2}(L^{2})}+\|\bfxi\|_{W^{2,2}(0,T)} \big)
\ea
\eeq{JB2}
and which, in addition, solves the original problem \eqref{4.7}. The proof of this convergence property is entirely analogous to that given in \cite[Lemma 3.4 and Section 4]{GS1}, to which we refer for the missing  details. Finally, the $T$-periodicity property of the pressure field is proved exactly as in \cite[Lemma 2]{GARMA}, and its proof will be omitted.  In order to complete the existence part of the lemma, we  recall some classical properties of  solutions to the Stokes problem:
\be
\ba{cc}\smallskip\left.\ba{ll}\medskip
\Delta\textbf{\textsf {w}}=\nabla {\sf p}+\textbf{\textsf{F}}\\
\Div\textbf{\textsf {w}}=0\ea\right\}\ \ \mbox{in $\Omega$}\\
\textbf{\textsf {w}}(x)=\textbf{\textsf {w}}_\star\,,\ \ x\in \partial\Omega\,.
\ea
\eeq{hm}
In particular,  we know that any distributional solution to \eqref{hm} satisfies the following estimate \cite[Lemma V.4.3]{Gab} 
\be
\|D^{2}\textbf{\textsf {w}}\|_2+\|\nabla {\sf p}\|_2\le C\,\big(\|\textbf{\textsf{F}}\|_{2}+\|\textbf{\textsf {w}}_\star\|_{3/2,2,\partial\Omega}+\|\textbf{\textsf {w}}\|_{2,\Omega_R}+\|{\sf p}\|_{2,\Omega_R}\big)\,,
\eeq{hm1}
with $C=C(\Omega,R)$.
Let $h\in L^{2}(\Omega_R)$ with $\int_{\Omega_R}h=0$, and let $\bfphi\in W^{1,2}_0(\Omega_R)$ be a solution to the problem $\Div\bfphi=h$ in $\Omega_R$, satisfying $\|\bfphi\|_{1,2}\le c_R\|h\|_{2}$. The existence of such a $\bfphi$ is well known \cite[Theorem III.3.1]{Gab}. Dot-multiplying both sides of \eqref{hm}$_1$ by $\bfphi$ and integrating by parts over $\Omega_R$, we get
$$
\langle \textbf{\textsf{F}},\bfphi\rangle+\langle\nabla\textbf{\textsf {w}},\nabla\bfphi\rangle=\langle{\sf p},\Div\bfphi\rangle=\langle {\sf p},h\rangle\,.
$$ 
From this relation,  the properties of $\bfphi$ and the arbitrariness of $h$, we deduce that ${\sf p}$, modified by a possible addition of a ($T$-periodic) function of time, must obeys the following inequality 
$$
\|{\sf p}\|_{2,\Omega_R}\le c_R\,\big(\|\textbf{\textsf{F}}\|_{2,\Omega_R}+\|\nabla\textbf{\textsf{w}}\|_{2,\Omega_R}\big)\,.
$$
As a result, \eqref{hm1}
furnishes
\be
\|D^{2}\textbf{\textsf {w}}\|_2+\|\nabla {\sf p}\|_2+\|{\sf p}\|_{2,\Omega_R}\le C\,\big(\|\textbf{\textsf{F}}\|_{2}+\|\textbf{\textsf {w}}_\star\|_{3/2,2,\partial\Omega}+\|\textbf{\textsf {w}}\|_{1,2,\Omega_R}\big)
\eeq{Li}
We next observe that, for each $t\in [0,T]$,  \eqref{4.7} can be put in the form \eqref{hm} with
$$
\textbf{\textsf{w}}\equiv \bfu\,,\ \  {\sf p}\equiv p\,,\ \ \textbf{\textsf{F}}\equiv \partial_t\bfu+\bfxi\cdot\nabla\bfu-\bff\,,\ \ \textbf{\textsf{w}}_\star\equiv \bfxi\,, 
$$
so that \eqref{Li} leads to
\be
\|D^{2}\bfu(t)\|_2+\|\nabla p(t)\|_2+\|p(t)\|_{2,\Omega_R}\le C_1\,\big(\|\bff(t)\|_{2}+|\bfxi(t)|+\|\partial_t\bfu(t)\|_{2}+\|\nabla\bfu(t)\|_{2}+\|\bfu(t)\|_{2,\Omega_R}\big)\,,
\eeq{Gy}
with $C_1=C_1(\Omega,R,\xi_0)$. If we combine \eqref{Gy} and use \eqref{JB2} we then show
\be
\|D^{2}\bfu\|_{L^\infty(L^2)}+\|\nabla p\|_{L^\infty(L^2)}+\|p(t)\|_{L^\infty(L^2(\Omega_R))}\le C\,\left(\|\bff\|_{W^{1,2}(L^2)}+\|\bfcalf\|_{L^{2}(L^{2})}+\|\bfxi\|_{W^{3,2}(0,T)} \right)\,.
\eeq{m1}
In view of \eqref{JB2} and \eqref{m1}, the proof of the existence property is thus completed. We shall now prove  uniqueness, namely,  that $\bfu\equiv\nabla p\equiv\0$ is the only $T$-periodic solution   in the class  \eqref{class} to the following system
\be\ba{cc}\smallskip\left.\ba{ll}\medskip
{\partial}_t\bfu-\bfxi(t)\cdot\nabla\bfu=\Delta\bfu-\nabla {p}\\
\Div\bfu=0\ea\right\}\ \ \mbox{in $\Omega\times (0,T)$}\\
\bfu(x,t)=\0\,,\ \ (x,t)\in \partial\Omega\times [0,T]\,.
\ea
\eeq{47}
To this end, we 
write
\be
\bfu=(\bfu-\bar{\bfu})+\bar{\bfu}:=\bfw+\bar{\bfu}\,,\ \ \bfxi=(\bfxi-\bar{\bfxi})+\bar{\bfxi}):=\bfchi+\bar{\bfxi}
. 
\eeq{471}
Since $\bar{\bfw}=0$, by Poincar\'e inequality,  Fubini's theorem and \eqref{class}, we deduce $\bfw\in L^{2}(L^2)$, so that, in particular,
\be
\bfw\in W^{1,2}(L^2)\cap L^2(W^{2,2})\,.
\eeq{2.15}
From classical embedding theorems (e.g. \cite[Theorem 2.1]{Sol}) and \eqref{2.15} we deduce
\be
\bfw\in L^\infty(L^2)\,.
\eeq{2.16}
Furthermore, from \eqref{47} it follows that $p$ obeys the following Neumann problem for a.a. $t\in [0,T]$
\be
\Delta p=0\ \ \mbox{in $\Omega$}\,,\ \ \pde p\bfn=-\curl(\psi\,\curl\bfu)\cdot\bfn\ \ \mbox{at $\partial\Omega$},
\eeq{2.17}
where $\psi$ is a smooth function of bounded support that is 1 in a neighborhood of $\partial\Omega$,
and we used the identity $\Delta\bfu=-\curl\curl\bfu.$ Employing well-known results on the Neumann problem \cite[Theorem III.3.2]{Gab} and the fact that $\bfu$ is in the class \eqref{class}, we get
$$
\|\nabla p\|_{L^2(L^q)}\le c\,\|\curl\curl\bfu\|_{L^2(L^q(K))}+\|\curl\bfu\|_{L^2(L^q(K))}\,,\ \ \mbox{all $q\in (1,2]$}\,,
$$
with $K=\supp(\psi)$. 
From this and Sobolev inequality, we may then modify $p$ by adding to it a suitable $T$-periodic function of time, in such a way that the redefined pressure field, that we continue to denote by $p$, satisfies
\be
p\in L^2(L^r)\,,\ \ \mbox{all $r\in(3/2,6]$}\,.
\eeq{2.20}
Let $\psi_R=\psi_R(x)$ be the function defined in \lemmref{1.1}.
We dot-multiply both sides of \eqref{47}$_1$ by $\psi_R\bfu$, and integrate by parts over $\Omega\times(0,T)$. Noticing that $\bfu\in L^2(L^2(\Omega\rho))$, all $\rho\ge R_*$, and using $T$-periodicity we thus show 
\be\ba{rl}\medskip
\Int0T\Int\Omega{}\psi_R\,|\nabla\bfu|^2&\!\!\!=-\half\Int0T\Int{\Omega^{\frac{R}{\sqrt2}}}{}\nabla\psi_R\cdot\bfxi(t)|\bfu|^2+\Int0T\Int{\Omega^{\frac{R}{\sqrt2}}}{} p\,\nabla\psi_R\cdot\bfu\\ &\!\!\!\!:= -\half I_{1R}+I_{2R}\,.
\ea\eeq{2.22}
From Schwarz inequality, the properties of $\psi_R$, and  \eqref{class} we get
$$
|I_{2R}|\le C_1\sup_{t\in [0,T]}\|\nabla\bfu(t)\|_2\int_0^T\|p(t)\|_{2,\Omega^{\frac{R}{\sqrt2}}}\,, 
$$
which, by \eqref{2.20}, entails
\be
\lim_{R\to\infty}|I_{2R}|=0\,.
\eeq{2.23}
Next, recalling that $\bar{\bfxi}=\lambda\,\bfe_1$, we may employ   \eqref{471} and Fubini's theorem to show
$$\ba{rl}\medskip
I_{1R}&\!=\Int{0}T\Int{\Omega_{\frac{R}{\sqrt2}}}{}\left[\lambda\,\pde{\psi_R}{x_1}\,|\bar{\bfu}|^2+\nabla\psi_R\cdot\bfxi(t)(|\bfw|^2+2\bar{\bfu}\cdot\bfw)\right]\\ &:=I_{1R}^1+I_{1R}^2\,,
\ea$$
where we have used the the fact that $\bar{\bfchi}=\0$. By H\"older inequality 
 and the summability properties of $\partial\psi_R/\partial x_1$ we show
$$
|I_{1R}^1|\le c\,\|\bfxi\|_{W^{1,2}(0,T)}\, \int_0^T\|\bar{\bfu}\|_{6,\Omega^{\frac{R}{\sqrt2}}}^2\,,
$$
which, 
in view of \eqref{class}, implies
\be
\lim_{R\to\infty}|I_{1R}^1|=0\,.
\eeq{2.24}
Finally, by using one more time Schwarz inequality and the properties of $\psi_R$, we infer
$$
|I_{2R}^2|\le 2\|\bfxi\|_{W^{1,2}(0,T)}\||\bfu|\,\nabla\psi_R\|_{2,\Omega^{\frac R{\sqrt2}}}\,\|\bfw\|_{L^\infty(L^2)}\le c\,\|\bfxi\|_{W^{1,2}(0,T)}\|\nabla\bfu\|_{2,\Omega^{\frac R{\sqrt2}}}\|\bfw\|_{L^\infty(L^2)}
\,,
$$
and so from the latter,  \eqref{2.16} and \eqref{class} we deduce
\be
\lim_{R\to\infty}|I_{2R}^2|=0\,.
\eeq{2.25}
Uniqueness then follows by letting $R\to\infty$ in \eqref{2.22} and using \eqref{2.23}--\eqref{2.25}. The lemma is completely proved.
\par\hfill$\square$\par
The following result provides, under further assumptions on $\bff$, the spatial asymptotic behavior of solutions determined in the previous lemma.
\Bl Let $(\bfu,p)$ be the solution to \eqref{4.7} constructed in \lemmref{4.2}. Then, if, in addition, $[\!]\bfcalf[\!]_{\infty,2,\lambda}<\infty$ it follows that  $[\!]\bfu[\!]_{\infty,1,\lambda}<\infty$, and, moreover,
$$
[\!]\bfu[\!]_{\infty,1,\lambda}\le C\,\big([\!]\bfcalf[\!]_{\infty,2,\lambda}+\|\bff\|_{W^{1,2}(L^2)}+\|\bfxi\|_{W^{2,2}(0,T)}\big)\,,
$$
where $C=C(\Omega,T,\xi_0)$, whenever $\|\bfxi\|_{W^{2,2}(0,T)}\in [0,\xi_0]$, for some $\xi_0>0$\,.
\EL{1}
{\em Proof.} Let $\psi$ be the ``cut-off" function introduced in \eqref{2.17}, and  let $\bfz$ be a solution to problem \eqref{Bog} with $f\equiv-\nabla\psi\cdot\bfu$. Since $\int_Kf=0$, where $K=\supp(f)$,  \lemmref{1.1_1} guarantees the existence of such a $\bfz$. Thus, setting
$$
\bfw:=\psi\,\bfu+\bfz\,,\ \ {\sf p}:=\psi\,p\,, \ \ \bfcalh=\psi\bfcalf
$$
from \eqref{4.7} we deduce that $(\bfw,{\sf p})$ is a $T$-periodic solution to the following problem
\be\left.\ba{ll}\medskip
{\partial}_t\bfw-\bfxi(t)\cdot\nabla\bfw=\Delta\bfw-\nabla {\sf p}+\Div\bfcalh+\bfg\\
\Div\bfw=0\ea\right\}\ \ \mbox{in ${\real^3}\times (0,T)$}
\,,
\eeq{1}
where
$$
\bfg:=-\partial_t\bfz+\bfxi(t)\cdot\nabla\bfz+\Delta\bfz-2\nabla\psi\cdot\nabla\bfu+p\,\nabla\psi-\bfxi(t)\cdot\nabla\psi\,\bfu\,.
$$
If we extend $\bfz$ to 0 outside its support, we infer that $\bfg$ is of bounded support. Also with the help of \lemmref{1.1_1} and \lemmref{4.2} we easily deduce
\be\ba{ll}\medskip
\Sup{t\ge 0}\|\bfg(t)\|_2\le c\,(\|\bff\|_{W^{1,2}(L^2)}+[\!]\bfcalf[\!]_{\infty,2,\lambda}+\|\bfxi\|_{W^{2,2}(0,T)})\,,\\
\Div\bfcalh(t)\in L^\infty(L^2)\,,
\ea
\eeq{g}
where we have used the obvious inequality $\|\bfcalf\|_{L^2(L^2)}\le c\,[\!]\bfcalf[\!]_{\infty,2,\lambda}$. We now introduce the new variable $\bfy$ defined by 
\be
\bfy=\bfx-\bfx_0(t)
\eeq{2}
where
\be
\bfx_0(t):=\int_0^t(\bfxi(s)-\bar{\bfxi})\,{\rm d}s\,.
\eeq{3}
Since $\bar{(\bfxi(t)-\bar{\bfxi})}=0$, one can show 
\be
\sup_{t\ge 0}|\bfx_0(t)|\le M
\eeq{GD0}
where
$$
M:=C\,\Frac{T^{\frac12}}{\pi}\big(\Int0T|\bfxi(t)-\bar{\bfxi}|^2\big)^{\frac12}
$$
and $C$ is numerical constant; see \cite{GARMA}. Thus, in particular,
\be
|\bfx|-M\le |\bfy|\le |\bfx|+M\,,
\eeq{GD}
Setting
\be\ba{ll}\medskip
\bfv(\bfy,t)=\bfw(\bfy+\bfx_0(t),t),\ \ 
{\sf P}(\bfy,t)={\sf p}(\bfy+\bfx_0(t),t), 
\\  
\bfcalg(\bfy,t)=\bfcalh(\bfy+\bfx_0(t),t)\,,\ \ \bfh=\bfg(\bfy+\bfx_0(t),t)
\ea
\eeq{4}
from \eqref{1} we easily deduce that $(\bfv,{\sf P})$ solves the following Cauchy problem
\be\ba{cc}\medskip\left.\ba{ll}\medskip
{\partial}_t\bfv-\lambda\,\partial_1\bfv=\Delta\bfv-\nabla {\sf P}+\Div\bfcalg+\bfh\\
\Div\bfv=0\ea\right\}\ \ \mbox{in ${\real^3}\times (0,\infty)$}
\,,\\
\bfv(x,0)=\bfw(x,0)\,.
\ea
\eeq{5}
We look for a solution to \eqref{5} of the form $(\bfv_1+\bfv_2, {\sf P}_1+{\sf P}_2)$ where
\be\ba{cc}\medskip\left.\ba{ll}\medskip
{\partial}_t\bfv_1-\lambda\,\partial_1\bfv_1=\Delta\bfv_1-\nabla {\sf P}_1+\Div\bfcalg+\bfh\\
\Div\bfv_1=0\ea\right\}\ \ \mbox{in ${\real^3}\times (0,\infty)$}
\,,\\
\bfv_1(x,0)=\0\,,
\ea
\eeq{6}
and
\be\ba{cc}\medskip\left.\ba{ll}\medskip
{\partial}_t\bfv_2-\lambda\,\partial_1\bfv_2=\Delta\bfv_2-\nabla {\sf P}_2\\
\Div\bfv_2=0\ea\right\}\ \ \mbox{in ${\real^3}\times (0,\infty)$}
\,,\\
\bfv_2(x,0)=\bfw(x,0)\,.
\ea
\eeq{7}
From \eqref{g} and \eqref{4}  we readily deduce
\be 
\ba{ll}\medskip
\Sup{t\ge 0}\,\|\bfh(t)\|_2\le c\,(\|\bff\|_{W^{1,2}(L^2)}+[\!]\bfcalf[\!]_{\infty,2,\lambda}+\|\bfxi\|_{W^{2,2}(0,T)})\,,\\
\Div\bfcalg(t)\in L^\infty(L^2)\,.
\ea
\eeq{g1}
Furthermore, by \eqref{GD0} and \eqref{GD} it follows that
\be\ba{rl}\medskip
(1+|x|)(1+2\lambda\,s(x))&\!\!\!\le (1+|y|+M)\big(1+2\lambda\,s(y)+2\lambda\,(M+x_{01}(t))\big)
\\
&\!\!\!\le C\,(1+|y|)\,\big(1+2\lambda\,s(y)\big)\,,
\ea
\eeq{LJ}
where here and in the rest of the proof $C$ denotes a constant depending, at most, on $\Omega$, $T$, and $\xi_0$. Likewise,
\be
(1+|y|)(1+2\lambda\,s(y))\le C\, (1+|y|+M)\big(1+2\lambda\,s(y)\big)
\,.
\eeq{VFM}
By \eqref{4} and the assumption on $\bfcalf$, \eqref{VFM} implies, in particular, $[\!]\bfcalg[\!]_{\infty,2,\lambda}<\infty$ and that
\be
[\!]\bfcalg[\!]_{\infty,2,\lambda}\le C\,[\!]\bfcalf[\!]_{\infty,2,\lambda}
\eeq{LJ1}
Thus, combining \eqref{g1}, \eqref{LJ1}and the assumption on $\bfcalf$ with \lemmref{VIII.4.4} we conclude that the Cauchy problem \eqref{6} has one (and only one) solution $(\bfv_1,{\sf P}_1)$ in the class (\ref{estimunn1}) for all $T>0$. Further,  we have $[\!]\bfv_1[\!]_{\infty,1,\lambda}<\infty$ with
\be
[\!]\bfv_1[\!]_{\infty,1,\lambda}\le C\,\big(\|\bff\|_{W^{1,2}(L^2)}+[\!]\bfcalf[\!]_{\infty,2,\lambda}+\|\bfxi\|_{W^{2,2}(0,T)}\big)
\eeq{FunEst}
Concerning \eqref{7}, a solution is given by
\be
v_{2i}(y,t)=\int_{\real^3}\Gamma_{i\ell}(y-z,s;\bar{\xi})\,w_{\ell}(z,0)\,{\rm d}z\,,
\eeq{14}
where $\bfGamma$ is the (time-dependent) Oseen fundamental tensor-solution to \eqref{7}$_{1,2}$  \cite[Theorem VIII.4.3]{Gab}.
Since $\bfw(x,0)\in L^6(\real^3)$, it follows that \cite[Theorem VIII.4.3]{Gab}
\be\ba{ll}\medskip
\bfv_2,\partial_t\bfv_2\, D^2\bfv_2\in L^r([\varepsilon,\tau]\times\real^3)\,,\ \ \mbox{all $\varepsilon\in (0,\tau)$, $\tau>0$, and $r\in[6,\infty]$}\,,\\ 
\|\bfv_2(t)\|_\infty\le C_1\,t^{-\frac14}\|\bfw(0)\|_6\,,\ \ \Sup{t\in(0,\infty)}\|\bfv_2(t)\|_6\le C_1\,\|\bfw(0)\|_6\,.
\ea 
\eeq{15}
In view of the regularity properties of $\bfu$ (and hence of $\bfw$) and those in (\ref{estimunn1}), \eqref{15} for $\bfv_i$, $i=1,2$, respectively, we may use the results proved in \cite[Lemma VIII.4.2]{Gab} to guarantee  $\bfw=\bfv_1+\bfv_2$. As a consequence, due to the $T$-periodicity of $\bfw$ and \eqref{4}$_1$, for  arbitrary positive integer $n$ and $t\in[0,T]$ we obtain 
\be\ba{rl}\medskip
|\bfw(x,t)|(1+|x|)(1+2\lambda\,s(x))&\!\!\!\!=|\bfv(y,t+nT)|(1+|x|)(1+2\lambda\,s(x))\\ &\!\!\!\!\le \big(|\bfv_1(y,t+nT)|+|\bfv_2(y,t+nT)|\big)(1+|x|)(1+2\lambda\,s(x)).\ea
\eeq{17}
Employing  \eqref{LJ}, \eqref{FunEst} and \eqref{15}$_2$ in this inequality we get
$$\ba{ll}\medskip
|\bfw(x,t)|(1+|x|)(1+2\lambda\,s(x))\\
\quad\quad
\le C\,\left[(1+|x|)(1+2\lambda\,s(x)) (t+nT)^{-\frac14}\|\bfw(0)\|_6+\|\bff\|_{W^{1,2}(L^2)}+[\!]\bfcalf[\!]_{\infty,2,\lambda}+\|\bfxi\|_{W^{2,2}(0,T)}\right]\,\ea
$$ 
so that,  by letting $n\to\infty$ and recalling that, uniformly in $t\ge 0$,  $\bfu(x,t)\equiv \bfw(x,t)$  for $|x|$ sufficiently large ($>\bar{R}$, say)  we deduce
\be
[\!]\bfu[\!]_{\infty,1,\lambda,\Omega^{\bar{R}}}\le C\,\big(\|\bff\|_{W^{1,2}(L^2)}+[\!]\bfcalf[\!]_{\infty,2,\lambda}+\|\bfxi\|_{W^{2,2}(0,T)}\big)\,. 
\eeq{vfn}
Since by classical embedding theorems we have
\be
\|\bfu\|_{L^\infty(L^\infty)}\le C\,\big(\|\bfu\|_{L^\infty(L^6)}+\|D^2\bfu\|_{L^\infty(L^2)}\big)\,,
\eeq{vfm}
the desired result then follows from \eqref{vfn}, \eqref{vfm} and \eqref{est}.
\par\hfill$\square$\par
The findings of \lemmref{4.2} and \lemmref{1} can be combined to arrive at the following one that represents the main achievement of this section.
\Bp Let 
$$\bff=\Div \bfcalf\in W^{1,2}(L^{2})\,, \ \ [\!]\bfcalf[\!]_{\infty,2,\lambda}<\infty\,,\ \ \bfxi\in W^{2,2}(0,T)
$$
be prescribed $T$-periodic functions with $\bar{\bfxi}=\lambda\,\bfe_1$, $\lambda\ge 0$. Then, there exists one and only one  $T$-periodic solution $(\bfu,p)$ to \eqref{4.7} such that
$$
[\!]\bfu[\!]_{\infty,1,\lambda}<\infty\,,\  \partial_t\bfu\in  L^{\infty}(W^{1,2})\cap L^2(D^{2,2})\,, \ \nabla\bfu\in L^\infty(W^{1,2})\,;\ \ p\in L^\infty(D^{1,2})\,.
$$
Furthermore, 
\be\ba{ll}\medskip
[\!]\bfu[\!]_{\infty,1,\lambda}+\|\partial_t\bfu\|_{L^{\infty}(W^{1,2})\cap L^2(D^{2,2})}+\|\nabla\bfu\|_{L^\infty(W^{1,2})} +\|\nabla p\|_{L^\infty(D^{1,2})}+\|p\|_{L^\infty(L^2(\Omega_R))}\\ 
\hspace*{1cm}\le C\,\big(\|\bff\|_{W^{1,2}(L^2)}+[\!]\bfcalf[\!]_{\infty,2,\lambda}+\|\bfxi\|_{W^{2,2}(0,T)} \big)
\ea
\eeq{est1}
where $C=C(\Omega,T,R,\xi_0)$, 
for any fixed $\xi_0$ such that $\|\bfxi\|_{W^{2,2}(0,T)}\le \xi_0$. 
\EP{1}
\setcounter{equation}{0}
\section{On the Unique Solvability of the Nonlinear Problem}
The main objective of this section is to study the properties of $T$-periodic solutions to the  full nonlinear problem \eqref{0.1}. 
This will be achieved by combining the results proved in \propref{1} with a classical contraction mapping argument. To this end, we introduce the Banach space
$$\ba{ll}\medskip
\mathscr S:=\big\{\mbox{$T$-periodic $\bfu$}: \Omega\times [0,T]\mapsto\real^3\,|\\
\qquad\qquad[\!]\bfu[\!]_{\infty,1,\lambda}<\infty\,,\  \partial_t\bfu\in  L^{\infty}(W^{1,2})\cap L^2(D^{2,2})\,, \ \nabla\bfu\in L^\infty(W^{1,2})\,;\ \Div\bfu=0\big\}\,, 
\ea$$
endowed with the norm
\be
\|\bfu\|_{\mathscr S}:=[\!]\bfu[\!]_{\infty,1,\lambda}+\|\partial_t\bfu\|_{L^{\infty}(W^{1,2})\cap L^2(D^{2,2})}+\|\nabla\bfu\|_{L^\infty(W^{1,2})}
\eeq{3.2}
\Bl Let $\textbf{\textsf{u}},\textbf{\textsf{w}}\in \mathscr S$. Then
$
\textbf{\textsf{u}}\cdot\nabla\textbf{\textsf{w}}\in W^{1,2}(L^2)$ and
$$
\|\textbf{\textsf{u}}\cdot\nabla\textbf{\textsf{w}}\|_{W^{1,2}(L^2)}\le c\,\|\textbf{\textsf{u}}\|_{\mathscr S}\|\textbf{\textsf{w}}\|_{\mathscr S}\,.
$$
\EL{3.1}
{\em Proof.} Clearly,
$$
\|\textbf{\textsf{u}}\cdot\nabla\textbf{\textsf{w}}\|_{L^{2}(L^2)}\le [\!]\textbf{\textsf{u}}[\!]_{\infty,1,\lambda}\,\|\nabla\textbf{\textsf{w}}\|_{L^\infty(L^2)}\le \|\textbf{\textsf{u}}\|_{\mathscr S}\|\textbf{\textsf{w}}\|_{\mathscr S}\,.
$$
Moreover, by using the embedding $L^4\subset W^{1,2}$ along with Schwarz inequality, we get
$$\ba{rl}\medskip
\|\partial_t\textbf{\textsf{u}}\cdot\nabla\textbf{\textsf{w}}\|_{L^{2}(L^2)}+\|\textbf{\textsf{u}}\cdot\nabla\partial_t\textbf{\textsf{w}}\|_{L^{2}(L^2)}&\!\!\!\le\|\partial_t\textbf{\textsf{u}}\|_{L^{4}(L^4)} \|\nabla\textbf{\textsf{w}}\|_{L^4(L^4)}+[\!]\textbf{\textsf{u}}[\!]_{\infty,1,\lambda}\, \|\nabla\partial_t\textbf{\textsf{w}}\|_{L^2(L^2)}\\ \medskip
&\!\!\! \le c\,\big(\|\partial_t\textbf{\textsf{u}}\|_{L^\infty(W^{1,2})}\|\nabla\textbf{\textsf{w}}\|_{L^\infty(W^{1,2})}+[\!]\textbf{\textsf{u}}[\!]_{\infty,1,\lambda}\|\partial_t\textbf{\textsf{w}}\|_{L^\infty(W^{1,2})}\big)\\
&\!\!\! \le c\,\|\textbf{\textsf{u}}\|_{\mathscr S}\|\textbf{\textsf{w}}\|_{\mathscr S}\,.
\ea
$$
The proof of the lemma is completed.
\par\hfill$\square$\par
We are now in a position to prove the main result of this paper.
\Bt
Let 
$$\bfb=\Div \bfcalb\in W^{1,2}(L^{2})\,, \ \ [\!]\bfcalb[\!]_{\infty,2,\lambda}<\infty\,,\ \ \bfxi\in W^{2,2}(0,T)
$$
be prescribed $T$-periodic functions with $\bar{\bfxi}=\lambda\,\bfe_1$, $\lambda\ge 0$.
Then, there exists  $\varepsilon_0=\varepsilon_0(\Omega,T)>0$ such that if 
$$
{\sf D}:=\|\bfb\|_{W^{1,2}(L^2)}+[\!]\bfcalb[\!]_{\infty,2,\lambda}+\|\bfxi\|_{W^{2,2}(0,T)}<\varepsilon_0\,,
$$
problem \eqref{0.1} has
one and only one  $T$-periodic solution $(\bfu,p)\in \mathscr S\times
 L^\infty(D^{1,2})$\,. Moreover, $\|\bfu\|_{\mathscr S}\le c\,{\sf D}$, for some $c=c(\Omega,T)$.
\ET{3.1}
{\em Proof.} We employ the contraction mapping theorem. To this end, define the map
$$
M:\textbf{\textsf{u}}\in\mathscr S\mapsto \bfu\in\mathscr S\,,
$$
with $\bfu$ solving the linear problem
\be\ba{cc}\smallskip\left.\ba{ll}\medskip
{\partial}_t\bfu-\bfxi(t)\cdot\nabla\bfu=\Delta\bfu-\nabla {p}+\textbf{\textsf{u}}\cdot\nabla \textbf{\textsf{u}}+\bfb\\
\Div\bfu=0\ea\right\}\ \ \mbox{in $\Omega\times (0,T)$}\\
\bfu(x,t)=\bfxi(t)\,,\ \ (x,t)\in \partial\Omega\times [0,T]\,,
\ea
\eeq{lin}
Set 
\be\textbf{\textsf{f}}:=\textbf{\textsf{u}}\cdot\nabla\textbf{\textsf{u}}=\Div(\textbf{\textsf{u}}\otimes \textbf{\textsf{u}}):=\Div\textbf{\textsf{F}}\,,
\eeq{C}
where we used the condition $\Div\textbf{\textsf{u}}=0$. In virtue of \lemmref{3.1},  by assumption, and by the obvious inequality
$$
[\!]\textbf{\textsf{F}}[\!]_{\infty,2,\lambda}\le c_1 [\!]\textbf{\textsf{u}}[\!]_{\infty,1,\lambda}^2\,,\ \ \textbf{\textsf{u}}\in\mathscr S\,,
$$
we infer that
$\textbf{\textsf{F}}$,  $\bfb$ and $\bfxi$ satisfy the assumption of \propref{1}. Therefore, by that proposition we conclude that the map $M$ is well defined and, in particular, that
\be
\|\bfu\|_{\mathscr S}\le c_2\left(\|\textbf{\textsf{u}}\|_{\mathscr S}^2+{\sf D}\right)\,,
\eeq{3.3}
with $c_2=c_2(\Omega,T,\xi_0)$. If we now take \be\|\textbf{\textsf{u}}\|_{\mathscr S}<\delta\,,\ \  \delta:=4c_2{\sf D}\,,\ \ {\sf D}<\frac{1}{16 c_2^2}\,,
\eeq{BD}
from \eqref{3.3} we deduce $\|\bfu\|_{\mathscr S}<\half\delta$.
Let $\textbf{\textsf{u}}_i\in\mathscr S$ $i=1,2$, and set 
$$
\textbf{\textsf{u}}:=\textbf{\textsf{u}}_1-\textbf{\textsf{u}}_2\,,\ \ \bfu:=M(\textbf{\textsf{u}}_1)-M(\textbf{\textsf{u}}_2)\,.
$$
From \eqref{lin} we then show
\be\ba{cc}\smallskip\left.\ba{ll}\medskip
{\partial}_t\bfu-\bfxi(t)\cdot\nabla\bfu=\Delta\bfu-\nabla {p}+\textbf{\textsf{u}}_1\cdot\nabla \textbf{\textsf{u}}+\textbf{\textsf{u}}\cdot\nabla\textbf{\textsf{u}}_2\\
\Div\bfu=0\ea\right\}\ \ \mbox{in $\Omega\times (0,T)$}\\
\bfu(x,t)=\0\,,\ \ (x,t)\in \partial\Omega\times [0,T]\,.
\ea
\eeq{line}
Proceeding as in the proof of \eqref{3.3} we can show
$$
\|\bfu\|_{\mathscr S}\le c_2\,\left(\|\textbf{\textsf{u}}_1\|_{\mathscr X}+\|\textbf{\textsf{u}}_2\|_{\mathscr S}\right)\|\textbf{\textsf{u}}\|_{\mathscr S}\,.
$$
As a result, if $\|\textbf{\textsf{u}}_i\|_{\mathscr S}<\delta$, $i=1,2$, from the previous inequality we infer
$$
\|\bfu\|_{\mathscr S}< 2c_2\delta\|\textbf{\textsf{u}}\|_{\mathscr X}\,,
$$
and since by \eqref{BD} $2c_2\delta<1/2$, we may conclude that $M$ is a contraction, which, along with \eqref{BD}, completes the proof of the theorem.
\par\hfill$\square$\par

\ed